\documentclass[12pt]{article}
\usepackage{amssymb}
\usepackage{amsmath}
\usepackage{amsbsy}
\usepackage{amsthm}
\usepackage{rotating}

\newcommand{\be}{\begin{equation}}
\newcommand{\ee}{\end{equation}}
\newcommand{\ba}{\begin{eqnarray}}
\newcommand{\ea}{\end{eqnarray}}
\newcommand{\ban}{\begin{eqnarray*}}
\newcommand{\ean}{\end{eqnarray*}}

\newtheorem{theo}{Theorem}[section]

\newtheorem{defi}[theo]{Definition}
\newtheorem{con}[theo]{Conjecture}

\usepackage{tikz}
\usetikzlibrary{matrix}
\usepackage{graphicx}

\begin{document}

\title{An inequality on global alliances for trees}
\author{Alexandria Yu}
\date{ }

\maketitle

{\noindent  Abstract:} In this paper, we prove an inequality on the cardinalities of the minimum size global defensive alliance and the minimum size global offensive alliance. A global defensive alliance is a dominating set such that when any point inside a selected group $S$ is chosen, at least half of the points in its neighborhood are also in the set $S$, including the selected point. A global offensive alliance is a dominating set such that if any point outside $S$ is selected, at least half of the points in its neighborhood, including the selected point, are in set $S$. Our result answers an open question in \cite{HA}.

\section{Introduction}

The concept of alliances can be used to show the agreements between groups, to predict how successful a new social networking site will be, and to illustrate the balance of power between different countries. \cite{EN} Alliances can also be used in determining business strategies among other things. \cite{EN} The idea of alliances was first introduced by Hedetniem-Hedetniem-Kristiansen \cite{HE} who introduced the concepts of global defensive and global offensive alliances. These sorts of alliances can be used to model classifications and online data flow. \cite{EN} Haynes-Hedetniem-Henning furthered this through studies of bipartite graphs and trees, contributing the discovery of their lower bounds \cite{HHH}. In his paper on bounds of global alliances of trees \cite{HA}, Harutyunyan proves that $|\gamma_o(T)-\gamma_a(T)|\leq\frac{n}{2}$ is true for all trees and proposed the following question:
Is it true that for any n-vertex tree $T, \gamma_o(T)\leq\gamma_a(T)+\frac{n}{6}$?

In this paper, we answer this question positively. Furthermore, we obtain an improved inequality, which is optimal.

The author would like to thank Sherry Gong for her helpful advice, interesting discussions, and encouragement.

\section{ Main result}

In this section, we state and prove the main theorem of this paper.

Let $G=(V, E)$ be a graph, where $V$ is the set of all vertices of $G$ and $E$ is the set of all edges of $G$.

\begin{defi}
 (1) For any subset $S$ of $V$, we define its boundary $\partial S$ to be the set of all vertices
  in $V-S$ which are adjacent to at least one vertex in $S$.
 (2) For any $v\in V$, we define $N[v]$ to be the neighborhood of $v$,
  so that $N[v]=v\cup\partial \{v\}$
  \end{defi}

\begin{defi} A dominating set is a set of points where for a subset $S$ of $V$, $S\cup\partial S=V$.

\end{defi}

\begin{defi} A set $S\subset V$ is called a defensive alliance if for every $v\in S$,
      $| N[v] \cap S| \geq |N[v]\cap (V-S)|$.  A defensive alliance $S$ is
       called a global defensive alliance if $S$ is also a dominating set.
\end{defi}

\begin{defi} A set $S\subset V$ is called an offensive alliance if for every $v\in \partial S$,
      $| N[v] \cap S| \geq |N[v]-S|$.  An offensive  alliance $S$ is
       called a global offensive  alliance if $S$ is also a dominating set.
\end{defi}

\begin{defi} The global defensive (offensive) alliance number of graph $G$ is the cardinality of a minimum
size global defensive (offensive) alliance in $G$, and is denoted $\gamma_a(G)$ ($\gamma_o(G))$. A minimum size global defensive
(offensive) alliance is called a $\gamma_a(G)$-set ($\gamma_o(G)$-set).
\end{defi}

\begin{con} 
As stated in \cite{HA}, if $T$ is a tree with $n$ vertices, then we have $$\gamma_a(T)+\frac{n}{6}\geq \gamma_o(T)$$.
\end{con}
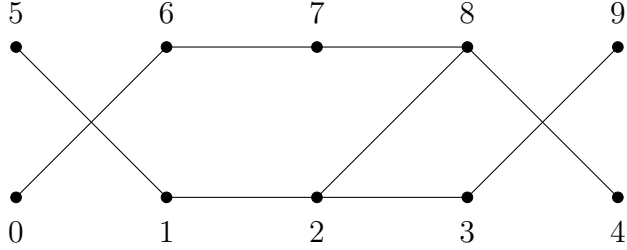
\begin{figure}
\centering
\begin{tikzpicture}

 \filldraw [black] (0,0) circle (2pt) node[below=2pt+3pt] {$0$}
                   (2,0) circle (2pt) node[below=2pt+3pt] {$1$}
                   (4,0) circle (2pt) node[below=2pt+3pt] {$2$}
                   (6,0) circle (2pt) node[below=2pt+3pt] {$3$}
                   (8,0) circle (2pt) node[below=2pt+3pt] {$4$}
                   (0,2) circle (2pt) node[above=2pt+3pt] {$5$}
                   (2,2) circle (2pt) node[above=2pt+3pt] {$6$}
                   (4,2) circle (2pt) node[above=2pt+3pt] {$7$}
                   (6,2) circle (2pt) node[above=2pt+3pt] {$8$}
                   (8,2) circle (2pt) node[above=2pt+3pt] {$9$}; 
\draw (0,0) -- (2,2) -- (4,2) -- (6,2) -- (8,0);
\draw (0,2) -- (2,0) -- (4,0) -- (6,0) -- (8,2);
\draw (4,0) -- (6,2);
\end{tikzpicture}
\caption{Example tree arranged in a bipartite graph.}
\label{example}
\end{figure}
\begin{theo} If $T$ is a tree with $n$ vertices, then we have $$\gamma_a(T)+\frac{n}{6}\geq \gamma_o(T).$$
\end{theo}
\proof
Step 1,
In this step, we prove that $\gamma_o(T)\leq \frac{n}{2}$. Given a tree $T$, the tree can be rearranged to form a bipartite graph. Select the side of the bipartite graph with the least number of vertices and call it $S$. This set $S$ is a global offensive alliance. This is because:
(1).	It is dominating: As the tree forms a bipartite graph, and is a tree, every non-selected vertex must be adjacent to a vertex in $S$.
(2).	It is offensive: As $T$ is bipartite, every point outside S will have exactly one neighbor also outside $S$ (as it counts itself as one neighbor) while, since $T$ is a tree, the vertex will be connected to at least one vertex in S. As S is the side of the bipartite graph with the least number of points, $S $ has at most $\frac{n}{2}$ points, where $n$ is the number of vertices in the tree.
Therefore, $|S|$  is less than or equal to $\frac{n}{2}$. That is, $\gamma_o(T)\leq \frac{n}{2}$

Step 2,
In this step, we will divide the situation into two cases to prove our final statement.

Case 1,

We will use a minimum size global defensive alliance to build an global offensive alliance by adding at most $\frac{n}{6}$ points to set $S$ of the offensive alliance.
 $$\gamma_a(T)=k\leq \frac{n}{3} +\frac{1}{2}.$$
Given a tree $T$, Let $k=|S|$ and let $E_S$ be the set of edges whose vertices lie in $S$.
Let $Y=V-S$.  Let $E_Y$ be the set of edges contained in $Y$.  let $E_B$ denote the set of all edges who has one vertex in
$S$ and another vertex in $Y$.
For any $x \in E_S$, by the defensive alliance condition, (the number of edges in $E_S$ with $x$ as one of its endpoints) + 1 $\geq$ the number of edges not in $E_S$  with $x$ as one of its endpoints.
Now summing over the vertices in $S$, we obtain $$2|E_S|+k\geq |E_B|$$
                   because each $e\in E_S$ is counted twice-- once for each of its vertices. 
As the tree is dominating, every vertex in $Y$ is connected to at least one other vertex in $S$. There are $n-k$ vertices in $Y$, so there are at least $n-k$ edges in  $E_B$.
It follows that $$|E_B|\geq n-k, so |E_S|\geq \frac{n}{2}-k.$$
Observe that the total number of edges in $T$ is $n-1$.
As a consequence, we obtain
$$|E_Y|=n-1-|E_S| -|E_B| \leq (n-1)-(\frac{n}{2}-k)-(n-k)=2k -\frac{n}{2}-1 ~~~~~~~~~~(\ast).$$

In order to build $S$ into a global offensive alliance, it suffices to "break" all of the edges contained in $Y$ by moving one vertex in every edge to $S$. Call the new set created by breaking the edges $S'$ and call its complement $Y'$. Since there are no edges remaining in $Y'$, all vertices located in $Y'$ will satisfy the condition for making $S'$ a global offensive alliance. By the above inequality $(\ast)$, we must add at most $|E_Y|\leq 2k -\frac{n}{2}-1$ number of vertices into $S$ in order to construct the global offensive alliance.

Since $ E_Y\leq 2k-\frac{n}{2} -1\leq \frac{n}{6},$ as $k\leq \frac{n}{3} +\frac{1}{2}$,
 the number of vertices we need to add to $S$ to get a global offensive alliance
 is less than or equal to $\frac{n}{6}$.

Case 2, 
 $$k>\frac{n}{3} +\frac{1}{2}.$$
 We have $$\gamma_a(T)=k> \frac{n}{3} +\frac{1}{2}.$$
It follows that $$\gamma_a(T)+\frac{n}{6} > \frac{n}{2}+\frac{1}{2}.$$

Recall that in Step 1, we proved that $$\gamma_o(T) \leq \frac{n}{2}.$$
As a consequence, we have $$\gamma_a(T)+\frac{n}{6} >  \gamma_o(T)+\frac{1}{2}.$$

And in both cases it is true that $$\gamma_o(T)\leq \gamma_a(T)+\frac{n}{6}.$$

Next, we can improve this inequality by adding the constant $C$ to the right side.
\begin{theo}If $T$ is a tree with $n$ vertices, then we have $$\gamma_a(T)+\frac{n}{6}\geq \gamma_o(T)+\frac{1}{3}.$$
\end{theo}
Now, we can divide the situation into two cases:
Case 1: If $k\leq \frac{n}{3}+C$, then $$E_Y\leq 2k-\frac{n}{2}-1\leq \frac{n}{6} +2C-1.$$
Hence $$\gamma_a(T)+\frac{n}{6}+2C-1 \geq   \gamma_o(T).$$
It follows that
$$\gamma_a(T)+\frac{n}{6} \geq \gamma_o(T)+1-2C.$$

Case 2: If $k> \frac{n}{3}+C$, then
$$\gamma_a(T)=k > \frac{n}{3}+C.$$
It follows that
$$\gamma_a(T)+\frac{n}{6} > \frac{n}{2} +C.$$
This implies that
$$\gamma_a(T)+\frac{n}{6} >  \gamma_o(T)+C.$$

We need to choose  a value of $C$ that maximizes the minimum between $C$ and $1-2C$.
As a result, two equations are produced:
\[ \begin{cases}
$y=1-2C$ \\
$y=C$
\end{cases} \]

\begin{figure}
\centering
\begin{tikzpicture}
  \draw[->] (-4,0) -- (4,0) node[right] {$c$};
  \draw[->] (0,-4) -- (0,4) node[above] {$y$};
  \draw[color=black]    plot[domain=-3:3] (\x,\x)             node[right] {$y=c$};
  \draw[color=black]    plot[domain=-0.5:2] (\x,{1-2*\x})             node[right] {$y=1-2c$};
  \filldraw (1/3, 1/3) circle (1pt);
  % \x r means to convert ?\x? from degrees to _r_adians:
%  \draw[color=orange] plot (\x,{0.05*exp(\x)}) node[right] {$f(x) = \frac{1}{20} \mathrm e^x$};
\end{tikzpicture}
 \caption{Graphed equations.}
\label{Graphed equations}
\end{figure}
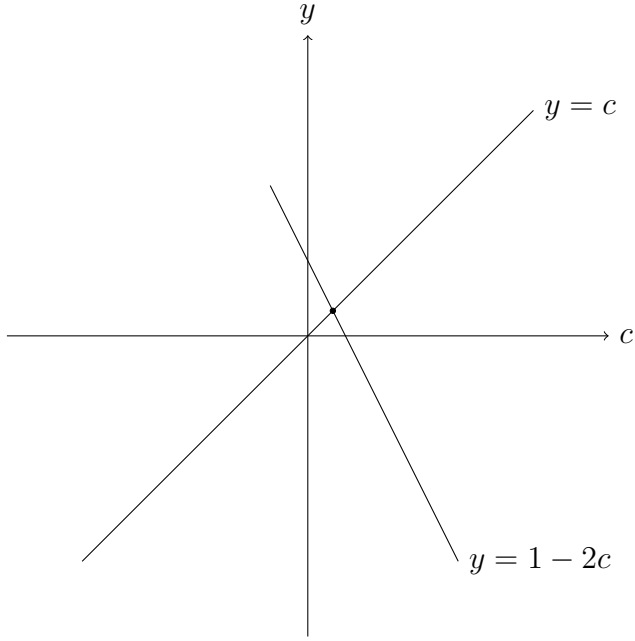
The intersection point is the point where this constant is optimal, thus by solving the equations, $\frac{1}{3}$ is the optimal value for $C$.

When $C=\frac{1}{3}$, we obtain $$\gamma_a(T)+\frac{n}{6} \geq  \gamma_o(T)+ \frac{1}{3}.$$
 Thus concludes the proof. \qed
\begin{figure}
\centering
\begin{tikzpicture}
 \filldraw [black] (1,0) circle (2pt) node[below=2pt+3pt] {$1$}
                   (4,0) circle (2pt) node[above=2pt+3pt] {$4$}
                   (5,0) circle (2pt) node[below=2pt+3pt] {$5$};
 \draw             (2,0) circle (2pt) node[below=2pt+3pt] {$2$}
                   (3,0) circle (2pt) node[below=2pt+3pt] {$3$}
                   (6,0) circle (2pt) node[above=2pt+3pt] {$8$}
                   (5,1) circle (2pt) node[above=2pt+3pt] {$7$}
                   (4,-1) circle (2pt) node[below=2pt+3pt]{$6$}; 
\draw  (1,0) -- (2,0) -- (3,0) -- (4,0) -- (5, 0) -- (6,0);
\draw  (4,0) -- (4,-1) (5,0) -- (5,1);
\end{tikzpicture}
\caption{Sharp-defensive.}
\label{sharpdefensive}
\end{figure}
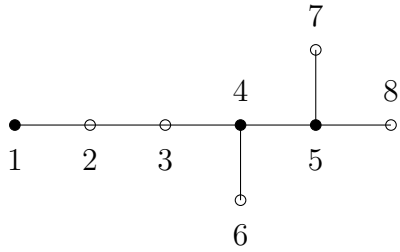

\begin{figure}
\centering
\begin{tikzpicture}
 \filldraw [black] (1,0) circle (2pt) node[below=2pt+3pt] {$1$}
                   (4,0) circle (2pt) node[above=2pt+3pt] {$4$}
                   (5,0) circle (2pt) node[below=2pt+3pt] {$5$}
                   (3,0) circle (2pt) node[below=2pt+3pt] {$3$};
 \draw             (2,0) circle (2pt) node[below=2pt+3pt] {$2$}
                   (6,0) circle (2pt) node[above=2pt+3pt] {$8$}
                   (5,1) circle (2pt) node[above=2pt+3pt] {$7$}
                   (4,-1) circle (2pt) node[below=2pt+3pt]{$6$}; 
\draw  (1,0) -- (2,0) -- (3,0) -- (4,0) -- (5, 0) -- (6,0);
\draw  (4,0) -- (4,-1) (5,0) -- (5,1);
\end{tikzpicture}
\caption{Sharp-offensive}
\label{sharp-offensive}
\end{figure}
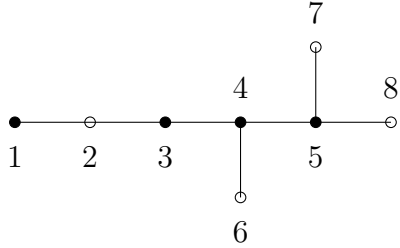

Note that our improvement on the inequality conjectured by Harutyunyan is sharp. For example, Figure \ref{sharpdefensive} shows a global defensive alliance that satisfies the equation $$\gamma_a(T)+\frac{n}{6} =\gamma_o(T)+ \frac{1}{3}$$ where the filled in vertices form the global defensive alliance and we add the vertex 3 to make a global offensive alliance, as shown in Figure \ref{sharp-offensive}.

Thus, $\gamma_a(T)=3$ $\gamma_o(T)=4$, and $\frac{n}{6}=\frac{8}{6}=\frac{4}{3}$. $\gamma_o(T)+\frac{n}{6}=3+\frac{4}{3}=4+\frac{1}{3}=\gamma_a(T)+\frac{1}{3}$, which satisfies our equation.

\noindent  University School of Nashville,

\noindent 2000 Edgehill Avenue, TN 37212,  USA.

\noindent e-mail: alexandriayu15@email.usn.org


\begin{thebibliography}{99}

\bibitem[HA]{HA} A. {\sc Harutyunyan}. \newblock {\em
 Some bounds on global allianes in trees.  } \newblock
 Discrete Applied Mathematics 161 (2013) 1739--1746.

\bibitem[HE]{HE} S.M. {\sc Hedetniemi}, S.T. {\sc Hedetniemi},  P. {\sc Krstiansen} \newblock {\em Alliances in graphs. } \newblock Journal of Combinatorial Mathematics and Combinatorial Computing 48 (2004) 105--177.

\bibitem[HHH]{HHH} T.W. {\sc Haynes},  S.M. {\sc Hedetniemi}, M.A. {\sc Henning}.   \newblock {\em Global defensive alliances in graphs. } \newblock
Eltronic Journal of Combinatorics, 10 (1) (2003) R47.

\bibitem[EN]{EN} A. R. (\sc Enciso). \newblock {\em
Alliances in graphs: Parameterized algorithms and on partitioning series- parallel graphs.} \newblock
 Dissertationas and Theses (2009)

\end{thebibliography}
\end{document}